\documentclass[conference]{ieeetran}
\usepackage{cite}
\usepackage{graphicx}
\usepackage{epstopdf}
\usepackage{amsmath}
\usepackage{cases}
\usepackage[caption=false,font=footnotesize]{subfig}
\usepackage{fixltx2e}
\usepackage{amssymb}
\usepackage{mathtools}
\usepackage{bm}
\usepackage{multirow}
\usepackage{color}
\usepackage{algorithm}
\usepackage{algorithmic}
\usepackage{amsthm}
\usepackage{mathrsfs}

\newtheorem{proposition}{Proposition}
\newtheorem{remark}{Remark}

\newcommand{\tabincell}[2]{\renewcommand\arraystretch{0.9}\begin{tabular}{@{}#1@{}}#2\end{tabular}}

\begin{document}

\title{Fast Monte Carlo Simulation of Dynamic Power Systems Under Continuous Random Disturbances}

\author{
  \IEEEauthorblockN{Yiwei~Qiu\textsuperscript{*}, Jin~Lin\textsuperscript{*}, Xiaoshuang~Chen\textsuperscript{*}, Feng~Liu\textsuperscript{*}, and Yonghua~Song\textsuperscript{\dag*}}
  \IEEEauthorblockA{
  \textsuperscript{*}State Key Laboratory of Control and Simulation of Power Systems and Generation Equipment, \\
  Department of Electrical Engineering, Tsinghua University, Beijing 100084, China\\
  \textsuperscript{\dag}Department of Electrical and Computer Engineering, The University of Macau, Macau 999078, China \\
  ywqiu@mail.tsinghua.edu.cn
  }
}

\maketitle

\begin{abstract}
    Continuous-time random disturbances from the renewable generation pose a significant impact on power system dynamic behavior. In evaluating this impact, the disturbances must be considered as continuous-time random processes instead of random variables that do not vary with time to ensure accuracy. Monte Carlo simulation (MCs) is a nonintrusive method to evaluate such impact that can be performed on commercial power system simulation software and is easy for power utilities to use, but is computationally cumbersome. Fast samplings methods such as Latin hypercube sampling (LHS) have been introduced to speed up sampling random variables, but yet cannot be applied to sample continuous disturbances. To overcome this limitation, this paper proposes a fast MCs method that enables the LHS to speed up sampling continuous disturbances, which is based on the It\^{o} process model of the disturbances and the approximation of the It\^{o} process by functions of independent normal random variables. A case study of the IEEE 39-Bus System shows that the proposed method is 47.6 and 6.7 times faster to converge compared to the traditional MCs in evaluating the expectation and variance of the system dynamic response.
\end{abstract}


\section{Introduction}
\label{sec:intro}

Continuous-time random disturbances due to the volatile renewable generation, e.g., wind and solar power, have a significant impact on power system dynamics. To evaluate this impact in system operation, the disturbances must be modeled with random processes instead of random variables that do not vary with time to ensure accuracy \cite{dong2012numerical,milano2013systematic,Chen2018Stochastic}. However, despite many studies on power system uncertainty evaluation considering random variables, e.g., probabilistic load flow (PLF), dynamic uncertainty evaluation methods for handling continuous disturbances are still inadequate.

Generally speaking, two types of methods, i.e., \emph{intrusive methods} \cite{Apostolopoulou2016assessment,ju2018analytical,shi2018analytical,li2019analytic,Wang2013the,Chen2018Stochastic}, and \emph{nonintrusive methods} \cite{dong2012numerical,milano2013systematic,wu2014stochastic} have been proposed to evaluate the impact of continuous random disturbances on dynamic power systems.

Intrusive methods use numerical computation program derived from highly specialized mathematical knowledge. In this sense, commercial power system simulation software such as \emph{PSS/E} cannot be employed. For example, in \cite{Wang2013the,Chen2018Stochastic} statistical information on power system dynamics under random disturbances is characterized by partial differential equations (PDEs) based on stochastic calculus. Clearly these PDEs cannot be solved by power system simulation software, making them impractical to be used by power utility companies.

In contrast, nonintrusive methods are based on commercial simulation software, which makes them easy for power utilities to use. However, to the authors' knowledge, currently Monte Carlo simulation (MCs) is the only available nonintrusive method for power system dynamic uncertainty evaluation in the presence of continuous random disturbances \cite{dong2012numerical,milano2013systematic,wu2014stochastic}, but the large number of sampling restricts its application.

To speed up MCs in dealing with random variables, for example in probabilistic load flow (PLF), Latin hypercube sampling (LHS) has been introduced, achieving much faster convergence than simple random sampling \cite{yu2009probabilistic,chen2012probabilistic}. Unfortunately, LHS cannot be directly used to sample random processes, therefore cannot be used to speed up the MCs of dynamic power systems under continuous random disturbances.

In this paper, to overcome this limitation, a fast MCs method for dynamic power systems under continuous random disturbances, which first enables LHS to deal with continuous random processes, is proposed.

Detailedly, the continuous random disturbances are approximated by functions of independent standard normal random variables based on the It\^{o} process model of the disturbances and the Karhunen-Lo\`{e}ve expansion (KLE) of Wiener processes that drive the It\^{o} process. Then, LHS is applied to sample the normal random variables, and disturbance signals are followingly reconstructed for simulation in commercial software. A case study in the IEEE 39-Bus System shows that the proposed method is $47.6$ times faster than the traditional MCs in evaluating expectation and variance, respectively.

\section{Problem Description}
\label{sec:problem}

\subsection{Uncertainty Assessment of Dynamics Power Systems Under Continuous Random Disturbances}

In evaluating the impact of continuous-time random disturbances from renewable generation on power system dynamics, the disturbances need to be modeled with random processes instead of static random variables that do not vary with time to guarantee accuracy \cite{dong2012numerical,milano2013systematic,Chen2018Stochastic}. Then, the dynamic response or performance index of the power system is determined by a function of the \emph{paths} of the random processes.

Specifically, a power system under continuous random disturbances can be depicted as differential-algebraic equations with the disturbances $\bm{\xi}_t$ as the parameters:
\begin{align}
      d{\bm{x}_t} &= \bm{f}(\bm{x}_t,\bm{y}_t;\bm{\xi}_t)dt, \label{eq:state} \\
      \bm{0} &= \bm{g} (\bm{x}_t,\bm{y}_t; \bm{\xi}_t),    \label{eq:algebra}
\end{align}
\noindent
where $\bm{x}_t$ and $\bm{y}_t$ are the states and algebraic variables; $t$ is time; and $\bm{f}(\cdot)$ and $\bm{g}(\cdot)$ are the state and algebraic equations.

Fixing certain initial condition, as the solution to (\ref{eq:state})--(\ref{eq:algebra}), system dynamic response such as the trajectory of rotor angle, or performance index such as CPS1 and CPS2 in automatic generation control (AGC) \cite{Chen2018Stochastic} can be described by an implicit function of the path of the excitation $\bm{\xi}_t$, denoted as
\begin{align}
    \omega = \omega \big(\left\{\bm{\xi}_\tau\right\}_{\tau\in[0,t]} \big), \label{eq:response}
\end{align}
\noindent
where $\omega(\cdot)$ is referred to as the \emph{random response function (RRF)} in this paper; and $\omega$ represents the value of $\omega(\cdot)$.

Clearly, as the paths of the disturbances randomly vary, the RRF varies accordingly. To provide assistance to system operation, statistical information on the RRF, e.g. expectation and variance, needs to be evaluated accurately and efficiently. This is the target of the uncertainty assessment.

\subsection{Motivation of the Proposed Method}

The renowned Monte Carlo simulation (MCs) can be used to evaluate the statistical information on the RRF \cite{dong2012numerical,wu2014stochastic}. In MCs, paths of the disturbances are repeatedly sampled and put into power system dynamic simulation software such as PSS/E to evaluate the RRF, and then statistical information is extracted. However, MCs can be computationally burdensome for practical application, since a large number of samplings are needed to achieve an acceptable precision.

To speed up sampling static random variables in MCs, for example in probabilistic load flow (PLF) problems, Latin hypercube sampling (LHS) is used \cite{chen2012probabilistic,wu2014stochastic}, which achieves a much faster convergence rate. However, currently LHS cannot be directly used to sample continuous-time random processes, which hinders it to be used to speed up MCs in the presence of continuous random disturbances.

\begin{figure}[t]
  \centering
  \includegraphics[width=3.5in]{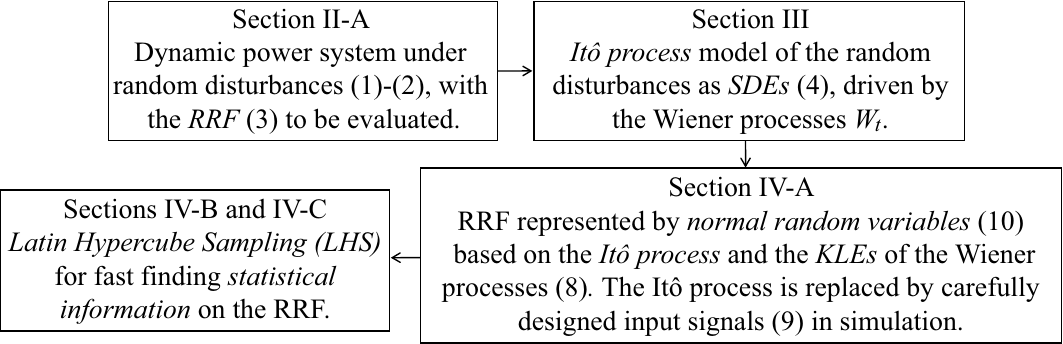}
  \caption{Brief framework of the proposed fast MCs method for power systems under continuous-time random disturbances.}
  \label{fig:frame}
\end{figure}

To overcome this problem, this paper uses the It\^{o} process to enable LHS to handle continuous random processes, and a fast MCs method for dynamic power systems under continuous random disturbances is proposed. The brief framework of the proposed method is shown in Fig. \ref{fig:frame}, and the detailed method are given in Sections \ref{sec:model} and \ref{sec:method}.

\section{It\^{o} Process Model of Continuous Random Disturbances}
\label{sec:model}

\subsection{Modeling Random Disturbances With the It\^{o} Process}
\label{sec:ito}

Our previous work \cite{Chen2018Stochastic} has shown that the volatile renewable generations can be modeled with an It\^{o} process, represented by the following stochastic differential equations (SDEs):
\begin{align}
  d\bm{\xi}_t & = \bm{\mu} \big( \bm{\xi}_t,t )dt + \bm{\sigma} \big( \bm{\xi}_t,t ) d{\bm{W}}_t,  \label{eq:ito}
\end{align}
\noindent
where $\bm{\xi}_t$ is the $m$-dimensional random disturbances; $\bm{W}_t = \left[ W_{1,t}, \ldots, W_{n,t} \right]^{\mathrm{T}}$ is the $n$-dimensional independent \emph{standard Wiener processes} (or \emph{Brownian motions}); $\bm{\mu}(\cdot): \mathbb{R}^m \times \mathbb{R}_+ \rightarrow \mathbb{R}^m$ and $\bm{\sigma}(\cdot): \mathbb{R}^m \times \mathbb{R}_+ \rightarrow \mathbb{R}^{m \times n}$ are called the \emph{drift} and \emph{diffusion terms}, respectively.

\begin{remark}
    By selecting different drift and diffusion terms, Gaussian or non-Gaussian random process $\bm{\xi}_t$ can be exactly characterized.
\end{remark}

For example, for one-dimensional It\^{o} processes subject to several typical probability distributions, the drift and diffusion terms are listed in Table \ref{tab:ito}. For other distributions, the corresponding It\^{o} process can also be constructed; for multidimensional random processes, their correlation is characterized by non-diagonal entries in the diffusion term; see \cite{Chen2018Stochastic} for details.

In addition, we present a method for identifying the It\^{o} process model from empirical data; see Appendix.


\begin{table}[tb]\scriptsize
  \renewcommand{\arraystretch}{1.45}
  \caption{Drift and Diffusion Terms of One-Dimensional It\^{o} Processes with Some Typical Probability Distributions}
  \label{tab:ito}
  \centering
  \begin{tabular}{cccc}
  \hline\hline
  Type          & Probability Density                                               & $\mu\left(\xi_t\right)$                & $\sigma^2\left(\xi_t\right)$        \\
  \hline
  Gaussian      & \tabincell{c}{$e^{{(\xi_t-a)^2}/2b}/\sqrt{2 \pi b}$}              & $-\left(\xi_t-a\right)$                & $2b$                  \\
  Beta          & \tabincell{c}{$\dfrac{\xi_t^{a-1}(1-\xi_t)^{b-1}}{B(a,b)}$}       & $-\left(\xi_t-\dfrac{a}{a+b}\right)$   & $\dfrac{2\xi_t\left(1-\xi_t\right)}{a+b}$        \\
  Gamma         & \tabincell{c}{$\dfrac{b^a}{\Gamma(a)}\xi_t^{a-1}e^{-b\xi_t}$}     & $-\left(\xi_t-a/b\right)$              & $2\xi_t/b$            \\
  Laplace       & \tabincell{c}{$e^{-|\xi_t-a|/b}/(2b)$}                            & $-\left(\xi_t-a\right)$                & $2b|\xi_t-a|+2b^2$    \\
  \hline\hline
  \end{tabular}
\end{table}

\subsection{Simulation of Continuous Random Disturbances}

In Monte Carlo simulation, paths of the random disturbances are created via a stochastic numerical integration scheme, for example the Maryuama-Euler (EM) scheme \cite{friz2010multidimensional}:
\begin{align}
    \bm{\xi}_{t + h} & = h \bm{\mu} \left( \bm{\xi}_t; \bm{q} \right) + \bm{\sigma} \big( \bm{\xi}_t; \bm{q} ) \sqrt{h}\bm{\zeta}, \label{eq:distito}
\end{align}
\noindent
where $h$ is the step length; $\bm{\zeta} \sim \mathcal{N}\left(\bm{0},\bm{I} \right)$ is sampled as a vector of independent normal random variables in each step.

The stochastic integration scheme (\ref{eq:distito}) also implies that time-domain discretization of the random disturbances is impractical for uncertainty evaluation since the number of the resulting random variables is proportional to the number of discretization steps, which causes the curse of dimensionality. However, discretizing the disturbances spectrally instead of in time domain can be a feasible solution. This is the basic idea of our proposed method illustrated later.

\section{The Proposed Fast Monte Carlo Simulation}
\label{sec:method}

\subsection{Spectral Representation of the Random Disturbances}
\label{sec:kl}

According to the Karhunen-Lo\`{e}ve theorem \cite{friz2010multidimensional}, a standard Wiener process $W_{i,t}$ can be spectrally decomposed into a series of independent standard normal random variables $\{\zeta_{i,j}\}_{j=1}^{\infty}$, i.e., the \emph{Karhunen-Lo\`{e}ve expansion (KLE)}, as
\begin{align}
  W_{i,T} = \int_0^T dW_{i,t} = \sum_{j=1}^{\infty} \zeta_{i,j} \int_0^T m_j(t) dt \label{eq:kl}
\end{align}
\noindent
where $\left\{m_j(t)\right\}_{j=1}^{\infty}$ are functions defined on interval $t\in[0,T]$:
\begin{align}
{m_j(t)=}
\begin{cases}
    \sqrt{{1}/{T}} &,\ j = 1\\
    \sqrt{{2}/{T}} \cos \left[ \dfrac{(j-1)\pi t}{T} \right] &,\ j \ge 2 \label{eq:klb2}
\end{cases}
\end{align}

For computation, the infinite series (\ref{eq:kl}) is truncated at a given order $K$. Then, taking the derivative of both sides of (\ref{eq:kl}) yields
\begin{align}
    dW_{i,t} \approx \sum_{j=1}^{K}  \zeta_{i,j} m_j(t). \label{eq:kle}
\end{align}

Substituting (\ref{eq:kle}) into the It\^{o} process (\ref{eq:ito}) yields the following ordinary differential equation with random coefficients:
\begin{align}
     \frac{d\bm{\xi}^*_t}{dt} (\bm{\zeta})& = \bm{\mu} \big( \bm{\xi}^*_t,t \big) + \bm{\sigma} \big( \bm{\xi}^*_t,t \big)  \sum_{j=1}^{K}  \bm{\zeta}_{j} m_j(t) \label{eq:eqkl}
\end{align}
\noindent
where $\bm{\zeta}_{j}$ is a vector of independent normal random variables with dimension $n$; and $\bm{\zeta}\triangleq \{\zeta_{i,j}\}_{1 \le i \le n,1 \le j \le K}$ is the vector of all independent normal random variables. For convenience, in the rest of this paper we reindex the entries of $\bm{\zeta}$ as $\{\zeta_i\}_{1\le i \le M}$ without ambiguity. $M=nK$ is the size of $\bm{\zeta}$.

By (\ref{eq:eqkl}), the continuous random disturbances are characterized by independent normal random variables. Then substitute (\ref{eq:eqkl}) into (\ref{eq:response}), the system dynamic response or performance index defined by the RRF $\omega(\cdot)$ can be approximated by a function (denoted as $\omega^*(\cdot)$) of the normal random variables:
\begin{align}
    \omega \approx \omega^*(\bm{\zeta}) = \omega\big(\left\{\bm{\xi}^*_\tau(\bm{\zeta})\right\}_{\tau\in[0,t]} \big). \label{eq:implict}
\end{align}
\noindent
This is ensured by the following proposition.
\begin{proposition}
\label{prop:2}
The approximate RRF $\omega^*(\bm{\zeta})$ in (\ref{eq:implict}) converges to the actual RRF $\omega(\bm{\zeta})^* \rightarrow \omega\big(\{\bm{\xi}_\tau\}_{\tau\in[0,t]}\big)$ , as $K \rightarrow \infty$.
\end{proposition}

Rigorous proof of this proposition can be obtained by combining the result of Section 15.5.3 of \cite{friz2010multidimensional} and the boundedness of the RRF (\ref{eq:response}). To save space we omit the detailed proof here.

This far, we have made the RRF be characterized by normal random variables as (\ref{eq:implict}). This allows LHS to be employed to speed up finding the statistical information on the RRF.

\subsection{Latin Hypercube Sampling of Random Variables in the KLE}
\label{sec:lhs}

\begin{figure}[tb]
  \centering
  \includegraphics[width=2.75in]{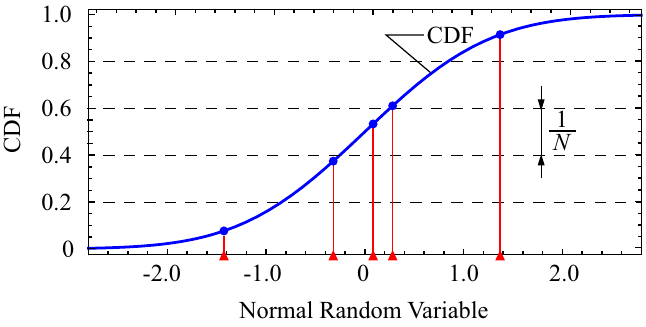}
  \caption{Sampling of a standard normal random variable in LHS.}
  \label{fig:lhs}
\end{figure}

The basic idea of Latin hypercube sampling (LHS) is to make the sampling point distribution close to the probability density function (PDF). Thus, fewer samplings are needed to achieve an acceptable precision in Monte Carlo simulation compared to simple sampling  \cite{yu2009probabilistic,chen2012probabilistic}. LHS has two steps:

\emph{Step 1: Random Sampling.} In our problem, the random variables are the coefficients  $\{\zeta_i\}_{1\le i \le M}$ in the KLE (\ref{eq:kl}).
Given sampling size $N$, for each random coefficient $\zeta_i$, evenly partition the range of the cumulative density function (CDF) into $N$ regions, and pick a random sample uniformly in each region. Then, sampling points are obtained by the inverse CDF of a standard normal distribution. This procedure is illustrated in Fig. \ref{fig:lhs}. The $N$ samplings of each random coefficient constitutes a row of an $M$-row primary sampling matrix.

\emph{Step 2: Permutation.} Random permutation on the primary sampling matrix is then performed. Permutation algorithms such as Cholesky decomposition \cite{yu2009probabilistic} can be employed to minimize the correlation between different columns. Due to space limit, here will not exposit the permutation algorithm. Interested readers are referred to the literature.

Once the two steps are finished, the $N$ samplings vectors are obtained as the $N$ columns of the permutated sampling matrix, denoted as $\{\hat{\bm{\zeta}}_k\}_{k=1}^{N}$. Then, paths of the random disturbances can be obtained using (\ref{eq:eqkl}) by replacing the random coefficients $\bm{\zeta}$ with the sampling vectors.

\begin{remark}
  Because the random variables $\bm{\zeta}$ in the KLE (\ref{eq:kl}) are independent, no additional transformation such as Nataf transformation is needed to make them independent. This makes our LHS much easier than that in correlated cases, for example in \cite{wu2014stochastic}.
\end{remark}

\subsection{Overall Computation Procedure}

The overall procedure of the proposed method is summarized as Algorithm \ref{alg:pipm}.

As clearly can be seen, the proposed method can be easily realized on commercial simulation software, making it very easy for power utilities to use in practical operation.

\begin{algorithm}[tb]
  \caption{Procedure of the Proposed Fast MCs Method}
  \label{alg:pipm}
  \begin{algorithmic}[1]
    \REQUIRE It\^{o} process model (\ref{eq:ito}) of the random disturbances, truncation order $K$ of the KLE (\ref{eq:kle}), sampling size $N$
    \STATE \label{step:1} create the Latin hypercube sampling set $\{\hat{\bm{\zeta}_k}\}_{k=1}^{N}$ following the procedure described in Section \ref{sec:lhs}.
    \STATE \label{step:2} for each sampling point $\hat{\bm{\zeta}}_k$, use (\ref{eq:eqkl}) to calculate the paths of corresponding disturbances
    \STATE \label{step:3} use power system simulation software such as PSS/E to find the RRF with respect to the paths obtained in Step \ref{step:2}
    \STATE extract statistical information on the RRF such as expectation and variance from the results obtained in Step \ref{step:3}
  \end{algorithmic}
\end{algorithm}

\section{Case Study}
\label{sec:case}

\subsection{Case Settings}
\label{sec:itomodel}

We compare the proposed fast dynamic uncertainty assessment method and the traditional Monte Carlo simulation method using the IEEE 39-bus system \cite{39sys}. The proposed method is coded in \emph{Python} with dynamic simulation performed on \emph{PSS/E}. Detailed models of GENROU generators, IEEET1 exciters, and TGOV1 governors are included.

\begin{figure}[tb]
  \centering
  \includegraphics[width=2.75in]{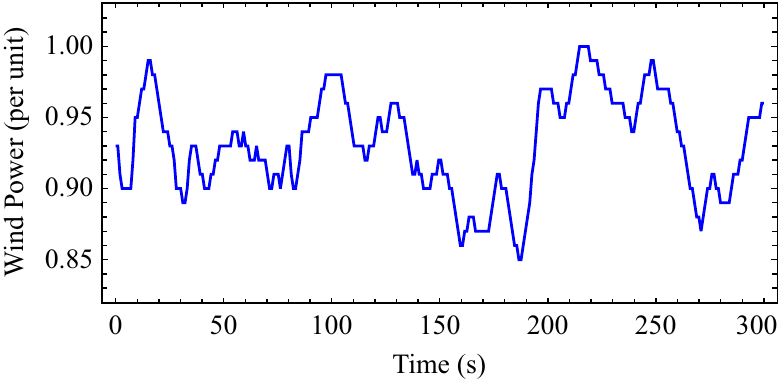}
  \caption{Wind power generation of an offshore wind farm over a $5$-minute interval, recorded by Ris\o\ DTU \cite{lin2012assessment}.}
  \label{fig:wind}
\end{figure}

A wind farm of rated power $3,000$ MW is connected to bus 15, which acts as a continuous random disturbance imposed on the system. The following It\^{o} process is used to model the per unit wind power, formulated as
\begin{align}
  dP_t = & \left[ 0.0535 - 0.0899 P_t + 0.0349 P_t^2 \right] dt  \nonumber \\
  & +  \left[-0.410 + 0.919 P_t - 0.505 P_t^2 \right] dW_t, \label{eq:windito}
\end{align}
\noindent
which is identified from the recorded data of an offshore wind farm \cite{lin2012assessment}, as shown in Fig. \ref{fig:wind}.

\begin{figure}[tb]
  \centering
  \includegraphics[width=3.42in]{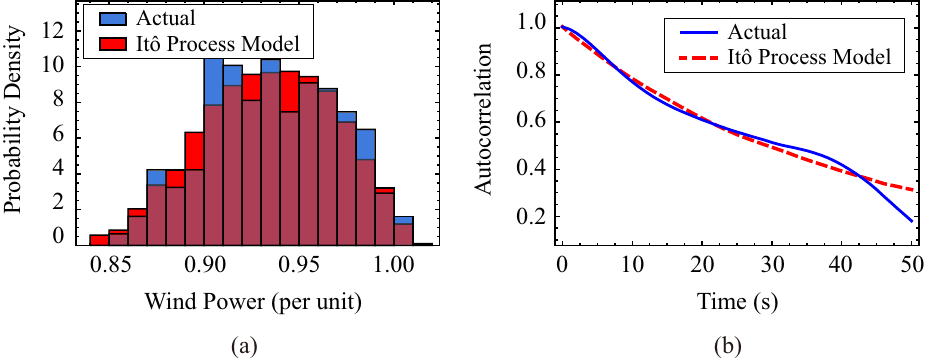}
  \caption{Probability density and autocorrelation of the actual wind power and the identified It\^{o} process model. (a) Probability density. (b) Autocorrelation.}
  \label{fig:wid}\vspace{-6pt}
\end{figure}

To validate the It\^{o} process model (\ref{eq:windito}), the probability distribution and autocorrelation of the recorded wind power data and simulated paths of the It\^{o} process model are compared in Fig. \ref{fig:wid}. Clearly, the result shows that non-Gaussian probability distribution and temporal correlation of wind power are precisely characterized by the identified It\^{o} process model.

\begin{figure}[tb]
  \centering
  \includegraphics[width=3.42in]{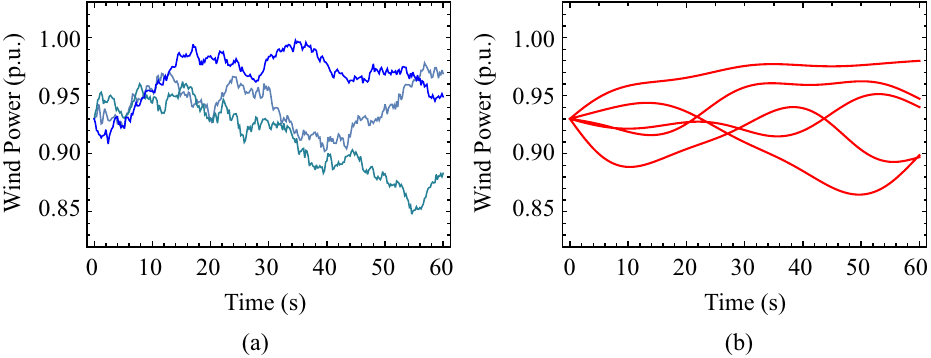}
  \caption{Paths of the random process of the wind power used for simulation. (a) Paths created by (\ref{eq:distito}) used in the traditional Monte Carlo simulation. (b) Paths created by (\ref{eq:eqkl}) used in the proposed method, with $K=6$.}
  \label{fig:cps}
\end{figure}

Next, we investigate the dynamic performance of frequency control under continuous wind power volatility. At $t=0$ s, the system is assumed at equilibrium, at $t=1.0$ s the generator at bus $30$ is tripped to simulate a fault. The RRF to be evaluated is defined as the root mean square (RMS) value of system frequency deviation within $60$ s.

\subsection{Comparison Between the Proposed Method and the Traditional Monte Carlo Simulation}
\label{sec:39}

\begin{figure}[tb]
  \centering
  \includegraphics[width=3.5in]{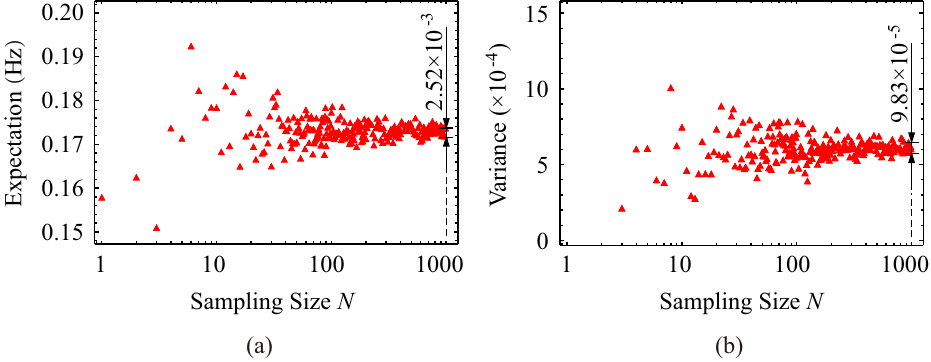}
  \caption{Expectation and variance of RMS frequency deviation obtained by the traditional MCs with different sampling size. (a) Expectation. (b) Variance.}
  \label{fig:mcs}
\end{figure}

\begin{figure}[tb]
  \centering
  \includegraphics[width=3.5in]{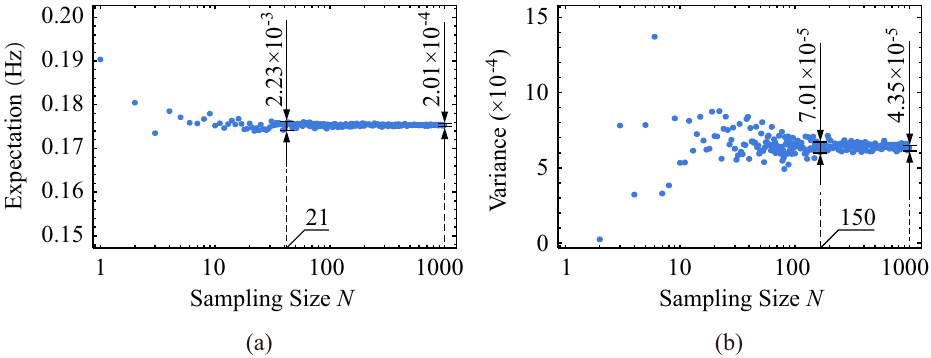}
  \caption{Expectation and variance of RMS frequency deviation obtained by the proposed method with different sampling size. (a) Expectation. (b) Variance. Clearly the convergence rate of the proposed method shown in this figure is much faster than the traditional MCs shown in Fig. \ref{fig:mcs}.}
  \label{fig:prop}
\end{figure}

In this section, the proposed method is compared to the traditional MCs to exhibit the improved efficiency. The two methods are respectively used to find statistical information, i.e., expectation and variance, on the RMS system frequency deviation. In MCs, paths of the It\^{o} process (\ref{eq:windito}) are sampled using the integration scheme (\ref{eq:distito}). In the proposed method, the order of KLE is set as $K = 6$, and the paths are created using (\ref{eq:eqkl}) with the random coefficients sampled via Latin hypercube sampling. For visualization, three sampled paths used in MCs are shown in Fig. \ref{fig:cps}(a), and five paths used in the proposed method are shown in Fig. \ref{fig:cps}(b).

With different sampling size $N$, the expectation and variance of the RMS frequency deviation obtained
by the traditional MCs are presented in Fig. \ref{fig:mcs}, and those obtained by the proposed method are shown in Fig. \ref{fig:prop}. Obviously, the proposed method achieves a much faster convergence than the traditional MCs either in evaluating the expectation or variance.

Quantitative comparison of the computational efficiency of the two methods is followingly performed. The maximal difference of the results with five successive sampling size is used to quantify the degree of convergence. With $N$ varies from $1$ to $1,000$, we find that in evaluating the expectation, the proposed method only need a sampling size of $N = 21$ to reach the degree of convergence of the traditional MCs with $N = 1,000$. This means the proposed method is $47.6$ times faster than the traditional MCs. As of variance, the proposed method also only needs $N = 150$ to reach the convergence degree of the traditional MCs with $N=1,000$,  about 6.7 times of efficiency compared to the traditional MCs. The related values of degree of convergence are labeled in Figs. \ref{fig:mcs} and \ref{fig:prop} respectively, and summarized in Table \ref{tab:conv} as well.

The above results effectively verify the improved efficiency of the proposed method compared to the existing MCs.

\begin{table}[tb]\scriptsize
  \renewcommand{\arraystretch}{1.4}
  \caption{Size of Sampling of the Proposed Method and the Traditional MCs for Convergence in Evaluating Expectation and Variance}
  \label{tab:conv}
  \centering
  \begin{tabular}{cccc}
  \hline \hline
                                & Method            & Sampling Size $N$ & Degree of Convergence \\ \hline
  \multirow{2}{*}{Expectation}  & Traditional MCs   & $1000$            & $2.52\times 10^{-3}$  \\
                                & The Proposed      & $21$              & $2.23\times 10^{-3}$  \\ \hline
  \multirow{2}{*}{Variance}     & Traditional MCs   & $1000$            & $9.83\times 10^{-5}$  \\
                                & The Proposed      & $150$             & $7.01\times 10^{-5}$  \\
  \hline \hline
  \end{tabular}
\end{table}

\section{Conclusions}

An uncertainty assessment method for dynamic power systems under continuous disturbances is proposed. The disturbances are approximated by functions of independent normal random variables, enabling  Latin hypercube sampling to be applied to speed up the Monte Carlo simulation.

Applying more advanced probabilistic analysis method on the proposed spectral approximation of continuous disturbance to achieve a faster and more precise assessment of variance and high-order moments is a promising future work.

\section*{Acknowledgement}

Financial supports from National Key R\&D Program of China (2018YFB0905200), National Natural Science Foundation of China (51907099, 51577096, 51677100, 51761135015), and China Postdoctoral Science Foundation are acknowledged.

\appendix[Method for Identifying the It\^{o} Process Model]
\label{sec:data}
\setcounter{equation}{0}
\renewcommand{\theequation}{\ref{sec:data}\arabic{equation}}

Suppose a set of recorded data with sampling interval $h$, denoted by $\{ \tilde{\bm{\xi}}_0, \tilde{\bm{\xi}}_h,\tilde{\bm{\xi}}_{2h},\ldots,\tilde{\bm{\xi}}_T \}$. Construct the drift and diffusion terms $\bm{\mu} \left(\bm{\xi}_t; \bm{q} \right)$ and $\bm{\sigma} \left(\bm{\xi}_t; \bm{q} \right)$ in (\ref{eq:ito}) as simple functions of $\bm{\xi}_t$ such as polynomials with parameters $\bm{q}$ to be identified, such that the likelihood of the following logarithmic conditional probability is maximized:
\begin{align}
    \max_{\bm{q}} L = \log \Pr \left[ \tilde{\bm{\xi}}_h,\tilde{\bm{\xi}}_{2h},\ldots,\tilde{\bm{\xi}}_T | \tilde{\bm{\xi}}_{0} \right].  \label{eq:likelihood}
\end{align}

By the independent incremental property of the It\^{o} process \cite{Pardoux2014Stochastic}, the conditional probability in (\ref{eq:likelihood}) can be reformed as
\begin{align}
    L = - \sum_{j=1}^{T/h} \log \Pr \left[ \tilde{\bm{\xi}}_{jh} | \tilde{\bm{\xi}}_{(j-1)h} \right]. \label{eq:like}
\end{align}

Considering the It\^{o} process (\ref{eq:ito}) in its discrete form (\ref{eq:distito}), and considering that the sampling interval $h$ is short, we obtain
\begin{align}
   \bm{\xi}_{t+h} \sim \mathcal{N} \left( \bm{\xi}_t + h \bm{\mu}\left( \bm{\xi}_t; \bm{q} \right), h\bm{\sigma}^\mathrm{T} \left( \bm{\xi}_t; \bm{q} \right)\bm{\sigma}\left( \bm{\xi}_t; \bm{q} \right) \right).
\end{align}

Therefore, the conditional probability in (\ref{eq:like}) is
\begin{align}
   \Pr & \Big[  \tilde{ \bm{\xi}}_{t+h}  | \tilde{\bm{\xi}}_{t}  \Big]   = \dfrac{1}{\sqrt{\left( 2 \pi h \right)^m \det\left( \tilde{\bm{\sigma}}_t\tilde{\bm{\sigma}}_t^\mathrm{T}\right) }}  \label{eq:reprob} \\
   \times   & \exp   \left\{ -\frac{1}{2}\left[  \Delta\tilde{\bm{\xi}}_{t} - h\tilde{\bm{\mu}}_t \right]^\mathrm{T} \left( \tilde{\bm{\sigma}}_t\tilde{\bm{\sigma}}_t^\mathrm{T}\right)^{-1} \left[ \Delta\tilde{\bm{\xi}}_{t} - h\tilde{\bm{\mu}}_t \right] \right\}, \nonumber
\end{align}
\noindent
where $\Delta \tilde{\bm{\xi}}_{t} = \tilde{\bm{\xi}}_{t+h} - \tilde{\bm{\xi}}_{t}$ represents the change in the recorded data over a sampling interval; $\tilde{\bm{\mu}}_t$ and $\tilde{\bm{\sigma}}_t$ represent $\bm{\mu}( \tilde{\bm{\xi}}_{t}; \bm{q})$ and $\bm{\sigma}( \tilde{\bm{\xi}}_{t}; \bm{q})$, respectively.

Substituting (\ref{eq:reprob}) into (\ref{eq:like}), letting $\tilde{\bm{D}}_{jh} =h \tilde{\bm{\sigma}}_{jh} \tilde{\bm{\sigma}}_{jh}^{\mathrm{T}}/2$, and neglecting the constant terms in the logarithmic function yields
\begin{align}
  \min_{\bm{q}} L'  = & \frac{1}{4} \sum_{j} \left[  \Delta\tilde{\bm{\xi}}_{jh} - h\tilde{\bm{\mu}}_{jh} \right]^\mathrm{T} \tilde{\bm{D}}_{jh}^{-1} \left[  \Delta\tilde{\bm{\xi}}_{jh} - h\tilde{\bm{\mu}}_{jh} \right] \nonumber \\
  + & \frac{1}{2}  \sum_{j} \log \det \left(\tilde{\bm{D}}_{jh}\right) \label{eq:ident}
\end{align}

Since (\ref{eq:ident}) is an unconstrained programming problem, common methods, such as gradient descent, can be used to find the optimal parameters for the It\^{o} process model.


\bibliographystyle{IEEEtran}
\bibliography{IEEEabrv,SAMC}

\end{document}